\documentclass[10pt]{article}
\usepackage{amssymb,amsmath,amscd,pb-diagram}
\usepackage[T1]{fontenc}
\usepackage{lscape}
\usepackage[all]{xy}

\usepackage{color}
\usepackage{epsfig}

\usepackage{pb-xy}

\newtheorem{theorem}{Theorem}[section]
\newtheorem{lemma}[theorem]{Lemma}
\newtheorem{proposition}[theorem]{Proposition}
\newtheorem{corollary}[theorem]{Corollary}

\newenvironment{definition}[1][Definition \addtocounter{theorem}{1}\thetheorem]{\begin{trivlist}\item [\hskip \labelsep {\bfseries #1}]}{\end{trivlist}}

\newenvironment{proof}[1][Proof]{\begin{trivlist}
\item[\hskip \labelsep {\textit{#1}}]}{\end{trivlist}}

\newenvironment{notation}[1][Notation]{\begin{trivlist}
\item[\hskip \labelsep {\bfseries #1}]}{\end{trivlist}}

\newcommand{\qed}{
\begin{flushright}
Q.E.D.
\end{flushright}
}
      
\newcommand{\Z}{\mathbb{Z}}

\newcommand{\R}{\mathbb{R}}

\newcommand{\U}{\mathbb{H}}

\newcommand{\pr}{PSL_2(\R)}

\newcommand{\srarrow}{\twoheadrightarrow}

\newcommand{\factor}[2]{\raisebox{0.1ex}{$#1$}\raisebox{-0.3ex}{/}\raisebox{-0.9ex}{$#2$}}

\newcounter{Lcount}

\title{A hyperbolic approach to $\exp_3(S^1)$}
\author{S. C. F. Rose}

\begin{document}

\maketitle

\begin{abstract}
In this paper we investigate a new geometric method of studying $exp_k(S^1)$, the set of all non-empty subsets of the circle of cardinality at most $k$. By considering the circle as the boundary of the hyperbolic plane we are able to use its group of isometries to determine explicitely the structure of its first few configuration spaces. We then study how these configuration spaces fit together in their union, $exp_3(S^1)$, to reprove an old theorem of Bott as well as to offer a new proof (following that of E. Shchepin) of the fact that the embedding $exp_1(S^1) \hookrightarrow exp_3(S^1)$ is the trefoil knot.
\end{abstract}

\section{Introduction and Preliminaries}

Since the mid 1900s, the study of finite subsets of topological spaces has cropped up in one guise or another---initially by Borsuk \cite{borsuk}, later corrected by Bott \cite{bott} and more recently by Handel \cite{handel}, Tuffley \cite{tuffley} and Mostovoy \cite{mostovoy}.

A few different definitions have been given in the literature (which are all easily seen to be equivalent), but the one that we will use for this paper is the following.

\begin{definition}
Let $X$ be a topological space, and let $\exp_k(X)$ denote the set of all non-empty subsets of $X$ of cardinality at most $k$. That is
\[
\exp_k(X) := \{ S \subseteq X \mid 0 < |S| \leq k\}
\]
There is a natural surjective map $\pi: X^k \srarrow \exp_k(X)$ which takes $(x_1, \ldots, x_k)$ to $\{x_1, \ldots, x_k\}$; we thus endow $\exp_k(X)$ with the quotient topology induced by this map. It is of course clear that under this topology, $\exp_1(X) \cong X$.
\end{definition}

There are a number of results that can be proven about this space and its relation to $X$. First of all, the map $X \mapsto \exp_k(X)$ is a homotopy functor from $U$, the category of topological spaces and continuous maps, to itself. More precisely, we have the following theorem.

\begin{theorem}
Let $X, Y$ be topological spaces, let $f, g : X \to Y$. If $F : X \times I \to Y$ is a homotopy from $f$ to $g$, then the map
\[
\overline{F} : \exp_k(X) \times I \to exp_k(Y) \quad \big(\{x_1, \ldots, x_k\}, t\big) \mapsto \{F(x_1, t), \ldots, F(x_k, t)\}
\]
is a homotopy from $\exp_k(f)$ to $\exp_k(g)$.
\end{theorem}

From which it immediately follows that

\begin{corollary}
The homotopy type of $\exp_k(X)$ depends only on that of $X$.
\end{corollary}

There are also natural inclusions $\exp_k(X) \hookrightarrow \exp_m(X)$ for $k \leq m$, taking the subset $\{x_1, \ldots, x_k\}$ to $\{x_1, \ldots, x_k\}$; thus we have
\[
X \cong \exp_1(X) \hookrightarrow \exp_2(X) \hookrightarrow \cdots \hookrightarrow \exp_k(X) \hookrightarrow \cdots
\]
with each $\exp_k(X)$---for $X$ Hausdorff---embedded as a closed subset of $\exp_m(X)$.

The proofs of these results (as well as many other topological and homotopy-theoretic results) are in \cite{handel}, which provides a good introduction to many of the basic properties of $\exp_k(X)$ (which he denotes as $Sub(X, k)$).

It should also be noted that much of this is in marked contrast to the notion of a configuration space, which is not functorial (unless we suitably restrict our category); nor are there nice embeddings as above.

The most famous result in the study of these spaces, however, comes from Bott's correction to Borsuk's initial paper. That is:

\begin{theorem}\label{themaintheorem}
The space $\exp_3(S^1)$ is homeomorphic to $S^3$, the three-sphere.
\end{theorem}

This has been proven so far in a number of different ways---Bott uses an elaborate ``cut-and-paste'' style argument (in effect showing that it the union of tori together with a calculation of its fundamental group to show that it is a simply connected lens space), Tuffley finds a decomposition of it as a $\Delta$-complex (in the sense of \cite{hatcher_book}) and shows that it is a simply connected Seifert-fibered space, while Mostovoy ties it in with a result of Quillen about lattices in the plane.

There are of course many other results about similar such spaces, primarily in \cite{tuffley} and \cite{sjerve_kallel}, most notably the following result.

\begin{theorem}\label{homotopy_type_of_exp_k}
The homotopy type of $\exp_k(S^1)$ is that of an odd dimensional sphere; specifically,
\[
\exp_{2k}(S^1) \simeq \exp_{2k-1}(S^1) \simeq S^{2k-1}
\]
\end{theorem}

This is proven in \cite{sjerve_kallel} using the notion of truncated product spaces---these can be considered to be, in a certain sense, the free $\factor{\Z}{2}$ set on points in our space. More specifically, these are given by
\[
TP_n(X) := \factor{SP_n(X)}{\sim}
\]
where $SP_n(X)$ is the $n$-th symmetric product of $X$ defined to be the collection of unordered $n$-tuples in $X$---that is, $SP_n(X) = \{x_1 + \cdots + x_n | x_i \in X\}$---and we say that $x_1 + \cdots + x_{k-2} + 2x \sim x_1 + \cdots x_{k-2}$ (where we use the convention of writing an unordered collection of points as a sum), topologized accordingly. The reason that this is useful is that $TP_n(S^1) \cong \R\mathbb{P}^n$ (see \cite{mostovoy_truncated}) and so if we then note the homeomorphism
\[
\factor{\exp_k(X)}{\exp_{k-1}(X)} \cong \factor{TP_k(X)}{TP_{k-2}(X)}
\]
this then becomes
\[
\factor{\exp_k(S^1)}{\exp_{k-1}(S^1)} \cong \factor{\R\mathbb{P}^k}{\R\mathbb{P}^{k-2}}
\]
We are then able to use the long exact sequence of homology groups to show that $\exp_k(S^1)$ is a homology sphere---since it is simply connected it is then homotopic to a sphere (by a standard argument).

Moreover, the fact that the $\exp$ functor is a homotopy functor has also been used in \cite{tuffley_graph} and \cite{tuffley_surface} to begin the determination of the homotopy type of $\exp_k(\Sigma)$ for $\Sigma$ a closed surface. The idea is that if you consider a punctured surface, then this is homotopic to a wedge of copies of $S^1$---and from the understanding of $\exp_k(S^1)$ Tuffley is then able to prove results about $\exp_k(\Gamma_\ell)$, where $\Gamma_\ell$ is a graph consisting of wedges of $\ell$ circles. From this, he uses a Mayer-Vietoris argument, covering a surface with $k + 1$ copies of $X_i = \Sigma \setminus \{p_i\}$ (for distinct points $p_i$), to determine information about the homology of the space $\exp_k(\Sigma)$ (this works since $\bigcap_{i=1}^{k+1} X_i = \Sigma \setminus \{p_1, \ldots, p_{k+1}\}$, and so these cover $\exp_k(\Sigma)$).

In this paper we intend to provide yet another proof Theorem \ref{themaintheorem}, beginning with certain properties of $\pr$. This method will generalize and provide methods of calculation not only of the homotopy type of higher dimensional analogues, but possibly even their homeomorphism type.

The main difference in this approach as compared to previous ones lies in its inherent geometric nature. This heavy emphasis on the underlying geometry of the circle (more generally of higher dimensional spheres) and its connection to hyperbolic space affords us a much more specific understanding of the topology of the resulting space. In contrast, previous methods (for higher dimensional analogues) were restricted to knowledge of the homology and homotopy groups of these spaces which permitted a knowledge of the homotopy type of these spaces and no more.

Specifically, we use the topology of the group of isometries of the hyperbolic plane $\U$---$\pr$---which we note is homeomorphic to $\U \times S^1$ via the map
\begin{equation}\label{homeo_type_of_pr}
T \mapsto \bigg(T(i), \frac{dT}{dz}(i)\bigg)
\end{equation}
Using this and the fact that $S^1 = \partial\U$, we then mod out certain subgroups which permits us to determine the homeomorphism type of $\exp_3(S^1)$ as claimed.

In this case, as in previous proofs (e.g. \cite{tuffley}), we take advantage of the machinery of Seifert fibred spaces. It is quite easy to see that $\exp_k(S^1)$ has an $S^1$ action on it given by
\[
\zeta \cdot \{x_1, \ldots, x_k\} = \{\zeta x_1, \ldots, \zeta x_k\}
\]
and so, given that $\exp_3(S^1)$ is a closed 3-manifold (lemma \ref{exp_is_a_manifold}), we can use well known classification theorems to determine the homeomorphism type of $\exp_3(S^1)$.

\begin{notation}
We will use the notation
\[
C_k(X) := \{S \subseteq X \mid |S| = k\}
\]
to denote the $k$-th unordered configuration space of $X$. It is of course clear that
\begin{equation}\label{exp_as_union}
\exp_k(X) = \bigcup_{j=1}^k C_j(X)
\end{equation}
The topic of this paper will largely be concerned with how this union behaves topologically.
\end{notation}



\section{The analysis of $\exp_3(S^1)$}

\subsection{The topology and geometry of $C_i(S^1)$}

The main structure of the proof of Theorem \ref{themaintheorem} is as follows. For any three points of $S^1 = \partial \U = \R \cup \infty$ there is an element of $\pr$  which takes the points to $\{0, 1, \infty\}$, preserving cyclic ordering. This element is not unique; post-composing it with any element permuting $\{0, 1, \infty\}$ yields another possible choice. If we let $\Gamma$ be the subgroup of $\pr$ generated by all elements permuting $\{0, 1, \infty\}$, then this subgroup acts on $\pr$ by left multiplication; if we quotient out by this action then we obtain the following result:

\begin{lemma} \label{group_theoretic_structure_of_c3}
The subgroup of $\pr$ which cyclically permutes $(0, 1, \infty)$ is simply the cyclic group $\Gamma = \factor{\Z}{3}$, generated by $\gamma = \tfrac{z - 1}{z}$; thus
\[
C_3(S^1) \cong \factor{\pr}{\Gamma}
\]
\end{lemma}

By similar reasoning we have that for any two points of $S^1$ there is an element of $\pr$ taking those points to $\{0, \infty\}$; if we then let $\Xi$ be the subgroup of $\pr$ which fixes the set $\{0, \infty\}$ we similarly have the following:

\begin{lemma} \label{group_theoretic_structure_of_c2}
The subgroup of $\pr$ which setwise fixes $\{0, \infty\}$ is generated by the elements $\tau = -\tfrac{1}{z}$ and $\sigma_\lambda = \lambda z$ (for $\lambda \in \R_{>0}$) subject to the relations
\begin{subequations}\label{c_2_relations}
\begin{align}
\tau\sigma_\lambda\tau   &= \sigma_\lambda^{-1} = \sigma_{\lambda^{-1}} \label{c_2_rel_commutativity}\\
\sigma_\lambda\sigma_\mu &= \sigma_{\lambda\mu} \\
\tau^2                   &= 1
\end{align}
\end{subequations}
which leaves us with
\[
C_2(S^1) \cong \factor{\pr}{\Xi}
\]
\end{lemma}

And of course, $C_1(S^1)$ is simply the circle $S^1$.


This tells us which subgroups we need to be modding out by; what we need to do next is study how these subgroups act on $\pr$ so that we may understand the topology of the resulting quotients.

\begin{proposition}\label{c2_mobius}
$C_2(S^1)$, described as the quotient $\factor{\pr}{\Xi}$ above, is homeomorphic to the open M\"obius band $\mathcal{M}$.
\end{proposition}

\begin{proof}
We have that $\Xi$ is generated by $\tau$ and the $\sigma_\lambda$, subject to the relations above. Note in particular that by \eqref{c_2_rel_commutativity} it suffices to consider the actions of $\sigma_\lambda$ and $\tau$ separately.

We consider the action of $\sigma_\lambda$ first. Now we have both
\begin{gather*}
(\sigma_\lambda \circ T)(i) = \lambda \big(T(i)\big) \\
\frac{d(\sigma_\lambda \circ T)}{dz}(i) = \lambda\frac{dT}{dz}(i)
\end{gather*}
and so the action of $\sigma_\lambda$ on $\pr \cong \U \times S^1$ is given by $\sigma_\lambda \cdot (z, \theta) = (\lambda z, \theta)$. If we then quotient out by the action of the subgroup of $\Xi$ generated by the $\sigma_\lambda$'s we find
\[
\factor{\pr}{\langle \sigma_\lambda \rangle} \cong (0, \pi) \times S^1
\]
where the angle $\phi \in (0, \pi)$ is simply $\arg(z)$.

Moving on to the action of $\tau$, we note that
\begin{gather*}
(\tau \circ T)(i) = -\frac{1}{T(i)} \\
\frac{d(\tau \circ T)}{dz}(i) = \frac{1}{T^2(i)}\frac{dT}{dz}(i)
\end{gather*}
and so since up to scaling (as an element of $(0, \pi)$) $-\frac{1}{T(i)} = \pi - \phi$, and similarly $\tfrac{1}{T^2(i)} = -2\phi$,  we find that
\[
\factor{\pr}{\Xi} \cong \factor{(0, \pi) \times S^1}{(\phi, \theta) \sim (\pi - \phi, \theta - 2\phi)}
\]
which can be pictorially represented as in figure \ref{why_its_a_mobius_band},  whence follows the conclusion.

\begin{figure}[ht]
\begin{center}
\input{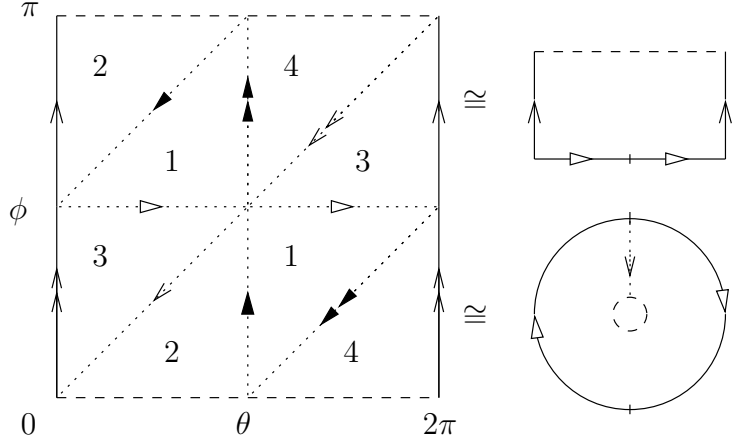}
\end{center}
\caption{The homeomorphism between $C_2(S^1)$ and $\mathcal{M}$}\label{why_its_a_mobius_band}
\end{figure}

\qed

\end{proof}

We next move on to the case of $C_3(S^1)$.

\begin{proposition}\label{homeo_type_of_c3}
$C_3(S^1)$, described as a the quotient $\factor{\pr}{\Gamma}$ above, is homeomorphic to the (open) model Seifert fibreing with twist $2\pi/3$.
\end{proposition}

\begin{proof}
Since $\gamma = 1 - \frac{1}{z}$, we have immediately that
\begin{gather}
(\gamma \circ T)(i) = 1 - \frac{1}{T(i)}\label{gamma_acts_part_1} \\
\frac{d(\gamma \circ T)}{dz}(i) = \frac{1}{T^2(i)}\frac{dT}{dz}(i)\label{gamma_acts_part_2}
\end{gather}
and so $\gamma$ does the same to the $S^1$ coordinate as $\tau$. Now, in the $\U$ coordinate, it turns out that $1 - \frac{1}{z}$ can be written as $R_1 \circ R_2$ where
\[
R_1(z) = 1 - \overline{z} \qquad R_2(z) = \frac{z}{|z|^2}
\]
and so it is the composition of two reflections, first about the circle $|z| = 1$ and then about the line $\Re z = \frac{1}{2}$---it is thus a rotation of the angle $\frac{4\pi}{3}$ about the point $e^{\pi i/3}$. With figure \ref{the_fundamental_domain4} as a fundamental domain, and with left and right edges identified according to equations \eqref{gamma_acts_part_1} and \eqref{gamma_acts_part_2}, we quickly see that this is a model Seifert-fibreing as claimed.

\begin{figure}[ht]
\begin{center}
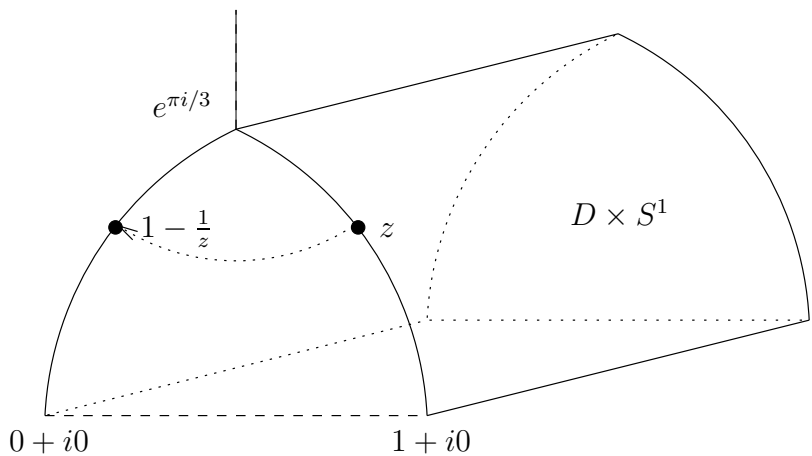
\end{center}
\caption{Choice of Fundamental Domain}\label{the_fundamental_domain4}
\end{figure}

\qed
\end{proof}

\subsection{Topology of the union $\exp_k(S^1)$}

At this point we have a solid grasp as to exactly what the three pieces $C_i(S^1), i = 1, 2, 3$ are. Remaining still is to describe the topology of their union---that is, what happens when points begin to coalesce.

We begin by looking at $\exp_2(S^1)$. First, in $\mathcal{M} \cong C_2(S^1)$, we find that a pair $\{p, r\}$ with $p > r$ maps to the element
\[
\xi(z) = \frac{z - p}{z - r}
\]
in $\pr$ (with a slight modification if $p, r = \infty$), and so using the homeomorphism from equation \eqref{homeo_type_of_pr} we see first that
\[
\xi(i) = \frac{i - p}{i - r} \qquad \frac{d\xi}{dz}(i) = \frac{p - r}{(i - r)^2}
\]
and so that as $p \to r$, $\xi(i) \to 1$, while $\frac{d\xi}{dz}(i) = \frac{p - r}{(i - r)^2}$, since we can ignore scaling, stays constant.

From this the following is immediate.

\begin{proposition}
The space $\exp_2(S^1) = C_1(S^1) \cup C_2(S^1)$ is homeomorphic to the closed M\"obius band $\mathcal{M}^*$.
\end{proposition}

\qed


We now move on to the situation with three points. First of all, our triple $\{p, q, r\}$ (with $p < q < r$, without loss of generality) maps to the element of $\pr$ given by
\[
\xi(z) = \left(\frac{q - r}{q - p}\right)\frac{z - p}{z - r}
\]
which, described in terms of equation \eqref{homeo_type_of_pr} yields
\begin{gather}
\xi(i) = \bigg(\frac{q - r}{q - p}\bigg)\frac{i - p}{i - r}\label{pqr_to_matrix} \\
\frac{d\xi}{dz}(i) = \bigg(\frac{q - r}{q - p}\bigg)\frac{p - r}{(i - r)^2}
\end{gather}
and so as in the case of two points, the position in the $S^1$ coordinate only depends on a single one of these points---in this case, $r$.



To aid in our description, we will assume that it is the point $q$ which tends towards either $p$ or $r$. That being the case, it should be noted that as $q$ varies, the quantity $\tfrac{q - r}{q - p}$ varies between $-\infty$ and $0$. As such (see equation \eqref{pqr_to_matrix}), by fixing both $p$ and $r$ we trace out a path in a particular $\U$ slice which happens to be a straight line from the point $0 + i0$ of slope $\tfrac{p - r}{1 + pr}$. Letting $q \to r$ (for simplicity's sake) we approach either $0 + i0$ or $1 + i0$, depending on the sign of $\tfrac{p - r}{1 + pr}$ (Consider our choice of fundamental domain). Then the M\"obius band is glued on to $C_3(S^1)$ entirely at the points $0 + i0$ and $1 + i0$. However, the point which we end up at on the M\"obius band depends on the slope of the path we take towards these points---we are ``blowing up'' these points---and so we find that these points are in fact separated by the set of lines passing through them, which are then glued to the corresponding point on the open M\"obius band.



But what about $C_1(S^1) = \partial(\mathcal{M}^*)$, and what of the edge running from $0 + i0$ to $1 + i0$? The first thing to note is that as the union $\exp_3(S^1)$ must be compact, any sequence $\{p_i, q_i, r_i\}$ which tends to that edge must necessarily converge to an element of $C_1(S^1)$ in the union. The only question remaining is what that element is.

Let $\{p_i, q_i, r_i\}$ be a sequence in $\exp_3(S^1)$ which converges (regarded as a point in the quotient $\factor{\U \times S^1}{\Gamma}$) to some point $P = (\lambda + i0, \theta)$ where $\lambda \in (0, 1)$, $\theta \in [0, \pi)$. Choose a neighbrouhood $N$ of $P$ such that $N$ is completely contained in the fundamental domain of figure \ref{the_fundamental_domain4}. As such, each $\{p_i, q_i, r_i\}$ is eventually contained in this neighbourhood, and thus in our particular choice of fundamental domain.

However, by our choice of fundamental domain there is an explicit ordering now on these points, and thus the $\U$ slice that any $\{p_i, q_i, r_i\} = \big[(p_i, q_i, r_i)\big]$ finds itself in is determined by a particular one of these points, say $r_i$. As this converges towards $(\lambda + i0, \theta)$, we must have that $r_i \to \theta$ (as well as $p_i, q_i$). That is, the point in $C_1(S^1)$ that we converge to is simply the point that all of $p_i, r_i, q_i$ are converging to. Thus in $\exp_3(S^1)$, the edge $(0 + i0, 1 + i0) \times \{\theta\}$ simply collapses to $\{*\}\times\{\theta\}$.

\subsection{$\exp_3(S^1)$ and the inclusion $\exp_1(S^1) \hookrightarrow \exp_3(S^1)$}

We will now prove (after a quick lemma) the main result of this paper, a stronger version of Theorem \ref{themaintheorem} (proven also in \cite{tuffley}).

\begin{lemma}\label{exp_is_a_manifold}
The space $\exp_3(S^1)$ is a compact $3$-manifold without boundary.
\end{lemma}

\begin{proof}
Compactness is immediate as $\exp_3(S^1)$ is a quotient of $(S^1)^3$. As for it being a manifold, the only place where this might fail is on $\exp_2(S^1) \subset \exp_3(S^1)$, so we simply need to verify that each point therein has a euclidean neighbourhood.

For points in $C_2(S^1)$ this is rather easy. A point in a neighbourhood of $\{p, q\}$ for $p \neq q$ will be in one of the configurations shown in figure \ref{nbhd_pq}, and so it is fairly easy to see that there is a neighbourhood of $\{p, q\}$ which is homeomorphic to two copies of $(-\epsilon, \epsilon) \times C_2\big((-\epsilon, \epsilon)\big)$ glued along their common boundary---which is then homeomorphic to $(-\epsilon, \epsilon)^3$ as required.

\begin{figure}[ht]
\begin{center}
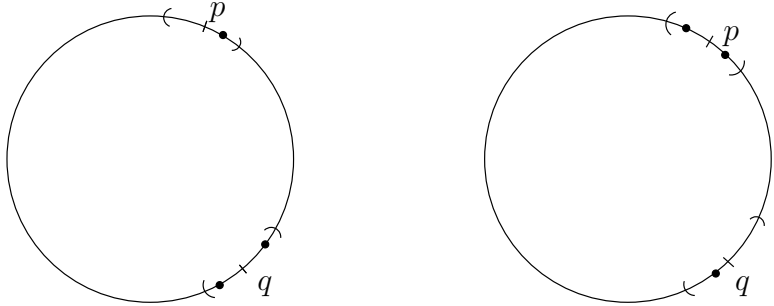
\end{center}
\caption{A neighbourhood of $\{p, q\}$}\label{nbhd_pq}
\end{figure}

For a point in $C_1(S^1)$, it is similarly easy to see that it has a neighbourhood homeomorphic to $\exp_3\big((-\epsilon, \epsilon)\big)$ which considered as a quotient of the space $X = \{(x, y, z) | -\epsilon < x \leq y \leq z < \epsilon\}$ (See figure \ref{exp_3_eps} with faces $A$, $B$ indentified) is also readily seen to be homeomorphic to a euclidean ball.

\begin{figure}[ht]
\begin{center}
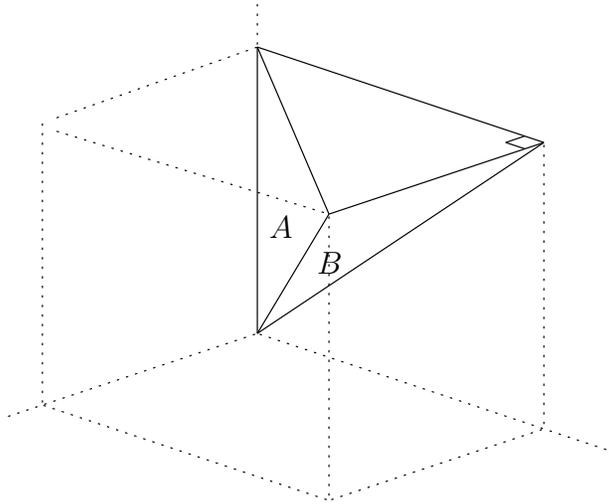
\end{center}
\caption{$\exp_3\big((-\epsilon, \epsilon)\big)$}\label{exp_3_eps}
\end{figure}

\qed
\end{proof}

\begin{theorem}
The space $\exp_3(S^1)$ is homeomorphic to $S^3$, and the inclusion $S^1 \cong \exp_1(S^1) \hookrightarrow \exp_3(S^1) \cong S^3$ is the trefoil knot.
\end{theorem}

\begin{proof}
The majority of this proof will rely on calculating the fundamental groups of $\exp_3(S^1)$ and $\exp_3(S^1) \setminus \exp_1(S^1)$, relying on classification theorems to show the above result. The main tool will be the Seifert-Van Kampen theorem.

From the covering shown in figure \ref{svk-domain_1}, the Seifert-Van Kampen theorem yields the following pushout of groups.

\[
\begin{diagram}
\node{\pi_1(A \cap B)}\arrow{e,t}{i_*}\arrow{s,l}{j_*} \node{\pi_1(B)}\arrow{s} \\
\node{\pi_1(A)}\arrow{e} \node{\pi_1\big(\exp_3(S^1)\big)}
\end{diagram}
\]
where $A$ deformation retracts onto a circle (and so $\pi_1(A) \cong \langle s \rangle$), and $B$ deformation retracts onto $\mathcal{M}^* \simeq S^1$ (hence $\pi_1(B) \cong \langle t \rangle$); thus we have that $\pi_1\big(\exp_3(S^1)\big) \cong \langle s, t \mid R \rangle$. It remains to determine what the relations $R$ are.

\begin{figure}[ht]
\begin{center}
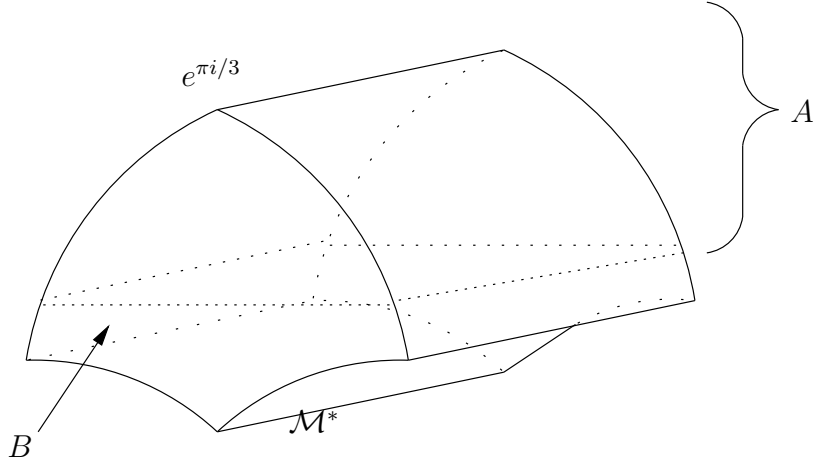
\end{center}
\caption{The covering used for the Seifert-Van Kampen theorem}\label{svk-domain_1}
\end{figure}

Now, $A \cap B$ is simply the ``boundary'' of $C_3(S^1)$; that is, up to homotopy it is simply a torus $T^2$. Thus its fundamental group is $\langle a \rangle \oplus \langle b \rangle$ where $a$ is the generator in the $S^1$ direction, and $b$ is the meridional generator.

We can explicitely describe this homotopy torus in terms of points on the underlying circle $S^1$ in the following manner. The longitudinal direction (ie the one corresponding to $a$ above) is given, as expected, simply by rotation of points along $S^1$; as we are avoiding the exceptional fibre, there is nothing unusual here and we need a full rotation at any given point to return to the starting point.

Now, we obtain the other generator $b$ (demonstrated in figure \ref{points_making_torus}) simply by rotating each point in a counter-clockwise manner to the next point along. It is easy to see that this commutes with $a$, and that together these two paths make a torus.


\begin{figure}[ht]
\begin{center}
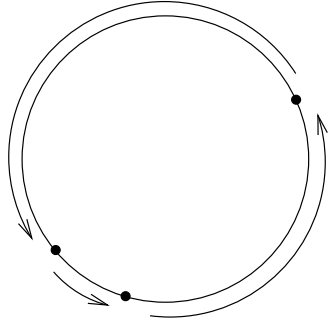
\end{center}
\caption{The generator $b$ of $\pi_1(A \cap B)$}\label{points_making_torus}
\end{figure}

So to determine the relations $R$, we first examine $i_*(a)$ and $j_*(a)$, the simplest of the two to deal with. Now, as the generator of $\pi_1(A)$ is the path along the exceptional fibre from $\theta=0$ to $\theta=\tfrac{\pi}{3}$, it is easy to see that $j_*(a) = s^3$. The exact same reasoning shows that $i_*(a) = t^2$, and so we have the relation that $s^3 = t^2$.

So for the meridional generator, $b$, we have the situation shown in figure \ref{mobius_seifert_path} (cf. figure \ref{why_its_a_mobius_band}) which can easily be seen to be homotopic to the generator $t$ of $\pi_1(B)$; thus it only remains to see what happens to $j_*(b)$ (shown in figure \ref{homotopy_j_of_b}) to fully understand what $\pi_1\big(\exp_3(S^1)\big)$ is. While it would be tempting to suggest that it simply collapses to a homotopically trivial path, this is not indeed the case.

\begin{figure}[ht]
\begin{center}
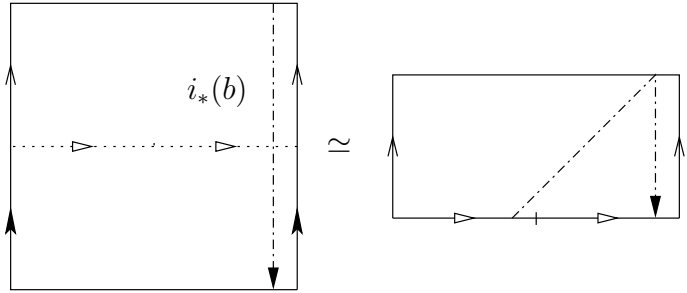
\end{center}
\caption{The homotopy class of $i_*(b)$}\label{mobius_seifert_path}
\end{figure}

\begin{figure}[ht]
\begin{center}
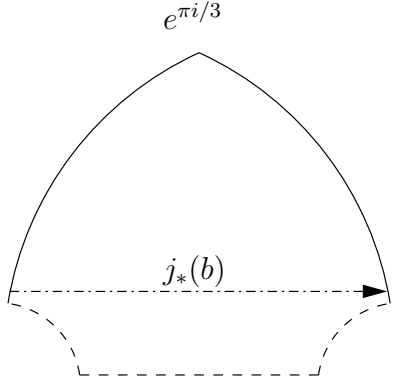
\end{center}
\caption{The homotopy class of $j_*(b)$}\label{homotopy_j_of_b}
\end{figure}

This is easily seen geometrically; if we perturb that diagram of figure \ref{points_making_torus} so that our three points are equally spaced about the circle, then it is easy to see that the description of the path $b$ is exactly the path which generates $\pi_1(A)$.

Combining all of this together, we find that
\[
\pi_1\big(\exp_3(S^1)\big) \cong \langle s, t \mid s^3 = t^2, s = t \rangle \cong 1
\]
and so $\exp_3(S^1)$ is simply connected. From \cite{seifert} it follows immediately that as a simply connected Seifert fibred space, this must be homeomorphic to $S^3$.

For brevity we will now define $X := \exp_3(S^1) \setminus \exp_1(S^1)$. Now, for the calculation of $\pi_1(X)$ the majority of the details above still hold through---the only difference is that what was labelled as $B$ above (now to be denoted $B'$) no longer deformation retracts onto $\mathcal{M}^*$.

Claim: $\pi_1(B') \cong \langle t, u \mid [t^2, u] = 1 \rangle$, where $t$ is the generator of $\pi_1(\mathcal{M})$ and $u$ is simply the image of the meridional generator in $B'$.

From this claim it follows that
\begin{align*}
\pi_1(X) &\cong \langle s, t, u \mid [t^2, u] = 1, s^3 = t^2, s = u \rangle \\
         &\cong \langle s, t \mid [t^2, s] = 1, s^3 = t^2 \rangle \\
         &\cong \langle s, t \mid s^3 = t^2 \rangle
\end{align*}

Now, to see that this implies that the inclusion $\exp_1(S^1) \hookrightarrow \exp_3(S^1)$ is the trefoil knot, we proceed as follows (This argument is due to E. Shchepin). The first thing is to note that the center of the M\"obius band---its exceptional fibre---is unkotted in $\exp_3(S^1)$. This is due to the fact that $\pi_1\big(C_3(S^1)\big) \cong \Z$. Thus if we consider a tubular neighbourhood of this subset we obtain a torus in $\exp_3(S^1)$---and the intersection of the boundary of this torus with the M\"obius band is thus a torus knot which is isotopic to $\exp_1(S^1)$. We can now use the fundamental group to say that it is a $(2, 3)$ torus knot, or a trefoil knot as claimed.

Proof of claim: Let us examine $B'$ a little more closely. Figure \ref{svk-domain_2} shows a slice of $B'$, separated into open sets $U, V$.

\begin{figure}[ht]
\begin{center}
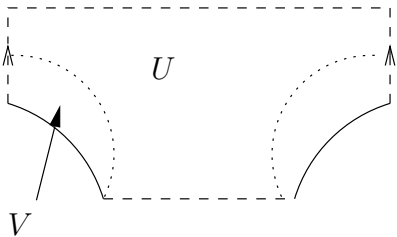
\end{center}
\caption{Covering of $B'$}\label{svk-domain_2}
\end{figure}

Now, $U \simeq T^2$, and $V \simeq \mathcal{M}$. Lastly, $U \cap V \simeq S^1$, and so we end up with the following pushout to calculate $\pi_1(B')$:

\[
\begin{diagram}
\node{\langle t \rangle}\arrow[2]{e,t}{t \mapsto a}\arrow{s,l}{t \mapsto c^2} \node[2]{\langle a \rangle \oplus \langle b \rangle}\arrow{s} \\
\node{\langle c \rangle}\arrow[2]{e} \node[2]{\pi_1(B')}
\end{diagram}
\]
from which it follows that
\begin{align*}
\pi_1(B') &\cong \langle a, b, c \mid [a, b] = 1, c^2 = a \rangle \\
          &\cong \langle b, c \mid [c^2, b] = 1 \rangle
\end{align*}
as claimed; Clearly, $c$ is the generator of $\pi_1(\mathcal{M})$, and $b$ is the meridional generator and so the full claim follows.
\qed
\end{proof}

\bibliography{../msc}

\end{document}